\documentclass[12pt]{amsart}

\usepackage{times,amsfonts,amsmath,amstext,amsbsy,amssymb,
  amsopn,amsthm,upref,eucal, amscd} 
\usepackage[T1]{fontenc}

\newtheorem{theorem}{Theorem}[section]
\newtheorem{lemma}[theorem]{Lemma}

\newtheorem{proposition}[theorem]{Proposition} 
\newtheorem{problem}{Problem}
\newtheorem{conjecture}[theorem]{Conjecture} 
\numberwithin{equation}{section}

\theoremstyle{definition}
\newtheorem{definition}[theorem]{Definition}

\newcommand{\field}[1]{\mathbb{#1}}
\newcommand{\R}{\field{R}} 
\newcommand{\N}{\field{N}}
\newcommand{\C}{\field{C}} 
\newcommand{\Z}{\field{Z}}
 
\newcommand{\T}{\field{T}} 
\newcommand{\Cal}{\mathcal}

\renewcommand{\>}{{\rangle}}
\newcommand{\pref}[1]{(\ref{#1})}

\begin{document}

\title[On  the Greenfield-Wallach and Katok conjectures]
{On  the Greenfield-Wallach and Katok \\ conjectures in dimension three}
\author{Giovanni Forni}

\address{Department of Mathematics, University of Toronto\\
40 St. George St., Toronto, Ontario,  CANADA M5S 2E4}

\email{forni@math.toronto.edu}

\address{Laboratoire de Math\'ematiques\\
        Universit\'e de Paris-Sud, B\^ atiment 425\\
       91405 Orsay Cedex, FRANCE}



\date{\today}

  


\maketitle

\section{Introduction} 

\noindent Many problems in dynamics can be reduced to the study of \emph{cohomological equations}
\cite{K01}, \cite{K03}. The simplest and most fundamental example of a cohomological equation for
a flow generated by a smooth vector field $X$ on a manifold $M$ is the linear partial differential equation 
\begin{equation}
\label{eq:CE}
Xu=f
\end{equation}
that is, the problem of finding a function $u$ on $M$ whose derivative along the flow is equal to
a given function $f$ on $M$. Roughly speaking, the $C^\infty$-cohomology of the flow is the space of non-trivial obstructions to the existence of a $C^\infty$ solution $u$ of \pref{eq:CE} for
any given function $f\in C^\infty(M)$. This notion is well-defined if the range of the Lie derivative
operator on the space $C^\infty(M)$ is closed. In this case the vector field is called \emph{$C^\infty$-stable}.

 In the 80's A.Katok \cite{H85}, \cite{K01}, \cite{K03} proposed the following conjecture. A vector field $X$ on a closed, connected orientable manifold is called \emph{cohomology free (CF) or rigid} if it is stable and the space of obstructions to the existence of solutions of the cohomological equation \pref{eq:CE} for smooth data is $1$-dimensional. A classical example of (CF) vector field, well-known from KAM theory, is given by constant  Diophantine vector fields on tori. Katok conjectured that these are the only examples up to smooth conjugacies.

 A related conjecture has been proposed earlier in 1973 by S. J. Greenfield and N. R. Wallach \cite{GW}. They introduced and studied \cite{GW}, \cite{GW2} the notion of a \emph{globally hypoelliptic (GH)} vector field  and conjectured that the only such vector fields (up to smooth conjugacies) are constant Diophantine vector fields on tori. A (GH) vector field $X$ is characterized by the property that if $XU$ is smooth for some distribution $U\in \Cal D'(M)$ then $U$ is smooth. This notion is modeled on the definition of a hypoelliptic differential operator in the theory of partial differential equations.

  In this paper we review recent progress, mainly due to W. Chen and M. Y. Chi \cite{CC00}, 
F. and J. Rodriguez-Hertz \cite{HH}, on the solution of these conjectures and derive a proof of both
the conjectures in dimension $3$. The argument reduces the (CF) conjecture in the general case to 
the contact case which can be finished by invoking the proof of the Weinstein conjecture recently announced by C. Taubes \cite{Taubes}. In fact, every (CF) flow is volume preserving and uniquely
ergodic, while according to the Weinstein conjecture every contact  flow on a closed, orientable
$3$-manifold has at least one periodic orbit.

In \cite{CC00} the authors propose a proof that every (GH) vector field on a torus is (CF). The argument, which is essentially correct and generalizes word for word to the general case, is presented below in \S \ref{sec:GH} . It follows that the notions of (GH) and (CF) vector fields and the related conjectures are equivalent.  Hence the Greenfield-Wallach conjecture in dimension $3$ is also proved. We are grateful to F. Rodriguez-Hertz who informed us of the results of \cite{CC00} in a personal communication. In this paper we prove:

\begin{theorem}  Let $M$ be a closed, connected, orientable $3$-manifold. If there exists a (CF) 
smooth vector field on $M$, then $M$ is diffeomorphic to a torus and $X$ is conjugate to a Diophantine 
vector field or $M$ is a rational homology sphere and $X$ is the Reeb vector field
of a smooth contact form. The latter case can be ruled out if the Weinstein conjecture holds. 
\end{theorem}

 The proof is based on the remarkable result of F. and J. Rodriguez-Hertz who proved
that if $M$ admits a (CF) vector field than $M$ fibers over a torus of dimension equal to the
first Betti number of $M$ \cite{HH}. Our argument consists in ruling out the existence of (CF) vector fields for all manifolds with Betti number strictly less than $3$. In case the Betti number is $2$ we prove that
every (CF) flow on $M$ is homogeneous, which contradicts a theorem of \cite{GW} which proves the
Greenfield-Wallach conjecture in the homogeneous case. In case the Betti number is equal to $1$, a
standard topological argument proves that any (CF) flow would have periodic orbits, a contradiction.
The case of vanishing Betti number is harder. A simple remark shows that any (CF) vector field is
tangent to a smooth plane field. In the integrable case, we are able to again prove that the flow is
homogeneous. In the non-integrable, contact, case, the Weinstein conjecture immediately implies the Greenfield-Wallach or the Katok conjecture. It would seem that a proof that there is no uniquely ergodic contact flow in dimension $3$ should be within reach of softer methods from the theory of dynamical
systems but we were so far unable to complete such an argument.

 We would like to acknowledge that partial proof of the results of this paper were obtained
independently by A. Kocsard in his Ph. D. thesis \cite{Ko}. In particular Kocsard independently 
proved the $3$-dimensional Katok conjecture for the cases of non-zero first Betti number. 
In the case of vanishing Betti number, after we informed him of our results, in particular of the 
existence of an invariant plane field, he produced an alternative proof in the integrable case. 
Finally, we would like to thank L. Flaminio for many discussions on the topics treated in the paper.

\section{Cohomology-free vector fields} 

\smallskip
\noindent Let $M$ be a closed, connected, orientable smooth manifold.
\begin{definition} 
\label{def:stable} (\cite{K01}, \cite{K03})
A smooth vector field $X$ on  $M$ is called \emph{$C^\infty$-stable} if the Lie derivative operator 
$\Cal L_X : C^\infty(M)  \to C^\infty(M)$ has closed range.
\end{definition}
\noindent If $X$ is stable, the $C^\infty$-cohomology of $X$ is well-defined and coincides with the space of all $X$-invariant distributions. Let $ \Cal D'(M)$ be the space of distributions on $M$ (in the sense of L. Schwartz), that is, the dual space of the Fr\'echet space $C^\infty(M)$.

\begin{definition}
A distribution $\Cal D \in \Cal D'(M)$  is called $X$-invariant if $X\Cal D=0$ in $\Cal D'(M)$. In other terms, the space $\Cal I_X(M)$ of $X$-invariant distributions is the kernel of the Lie derivative operator $\Cal L_X :\Cal D'(M) \to \Cal D'(M) $. 
\end{definition}
\noindent By definition a distribution is $X$-invariant if and only if it vanishes on the range of the operator $\Cal L_X : C^\infty(M) \to C^\infty(M)$. It follows from the Hahn-Banach theorem that if $X$ is $C^\infty$-stable,  the cohomological equation \pref{eq:CE} has a solution $u\in C^\infty(M)$ if and only if $\Cal D(f)=0$ for all $\Cal D\in \Cal I_X(M)$. 

\noindent Similar notions of stability and cohomology of a vector field can be introduced for different regularity classes \cite{K01}, \cite{K03}. For instance, we can say that $X$ is $(r,s)$-stable if the set 
$\{ Xu \in C^s(M) \,\vert \, u\in C^r(M)\}$ is closed in $C^s(M)$. The $(r,s)$-cohomology of $X$ is then the set of obstructions to the existence of a solution $u\in C^r(M)$ for a given $f\in C^s(M)$, that is, the subspace of $X$-invariant distributions which belong to the dual space $C^s(M)^\ast$. In this paper we will consider only the case $r=s=\infty$. 

\noindent  It is clear that all Borel probability measures invariant under the flow $\{\phi_t^X\}$ generated by $X$ are $X$-invariant distributions via integration. By the Krylov-Bogoliubov's theorem if $M$ is compact there exists at least one invariant probability measure for any continuous flow on $M$. It follows that the range of the operator $\Cal L_X$ on $C^\infty(M)$ has codimension at least $1$. 

\begin{definition} (\cite{K01}, \cite{K03})
A smooth vector field $X$ on $M$ is \emph{ $C^\infty$-cohomology free  (CF)} or \emph{$C^\infty$-rigid}, if for all $f\in C^\infty(M)$ there exists a constant $c(f)\in \C$ and $u\in C^\infty(M)$ such that
  $$ Xu = f -c(f). $$
\end{definition}
\noindent It is immediate to verify that the properties of stability and rigidity are invariant under 
$C^\infty$ conjugacies. The fundamental dynamical properties of (CF) vector fields are easily 
proved.

\begin{proposition}  \label{prop:CFdyn} [\cite{K01} p. 21]  Let $X$ be a (CF) vector field. Then 
the flow $\{\phi_t^X\}$ is conservative, that is, there exists a smooth $\{\phi_t^X\}$-invariant volume 
form $\omega$ on $M$. The space $\Cal I_X(M)$ of $X$-invariant distributions is one-dimensional 
and equal to $\C \omega$. In particular the flow $\{\phi_t^X\}$ is uniquely ergodic and minimal (strictly
ergodic).
\end{proposition}
\begin{proof}  It is immediate to prove that $\Cal I_X(M)$ is one-dimensional, hence in particular
the flow $\{\phi_t^X\}$ is uniquely ergodic. Let $w$ be any smooth volume form and let $f\in C^\infty(M)$ be the function such that $\Cal L_X w = f w$. Since $X$ is (CF), there exists $c_f \in \C$ and $u\in C^\infty(M)$ such that $Xu=f- c_f$. Let $\omega:= e^{-u} w$. Then 
$$
\Cal L_X \omega= (f-Xu)\omega  = c_f \omega \,.
$$
This implies that  $c_f =0$ since
$$ 
c_f \text{ \rm vol }_\omega (M) = \int_M \Cal L_X \omega = \int_M d (\imath_X \omega) = 0\,.
$$
It follows that the volume form $\omega$ is $\{\phi_t^X\}$-invariant, hence it coincides,
up to normalization, with the unique invariant probability measure. 
\end{proof}

\noindent If $M=\T^n$ is the $n$-dimensional torus, it is a simple but fundamental result that all
\emph{Diophantine} constant  flows are (CF), while all ergodic \emph{Liouvillean} constant  flows
are not $C^\infty$-stable (see \cite{K01} p. 19).  We recall that a constant  vector field $X=(\alpha_1, \dots, \alpha_n)$ on $\T^n$ is called Diophantine if there exist constants $\gamma >0$ and 
$C>0$ such that
$$
\vert \sum_{i=1}^n k_i \alpha_i \vert \geq   \frac {C }{ \Vert k \Vert ^\gamma} \quad
\text{ \rm for all } k=(k_1, \dots, k_n) \in \Z^n\setminus\{0\}\,.
$$
An ergodic constant  vector field on $\T^n$ that is not Diophantine is called \emph{Liouvillean}. The cohomological equation for constant vector fields on $\T^n$ can be analyzed by means of the standard Fourier series expansions. A stronger result which can be derived by adapting to flows methods developed for maps by Luz and dos Santos \cite{LdS}  is that every (CF) vector field on a torus is smoothly conjugate to a constant Diophantine vector field (see \cite{Ko}, \S 2.2). A more general 
result which holds for any closed, connected orientable manifold has been proved recently by F.  
and J. Rodriguez-Hertz \cite{HH}. Their work will be outlined below in \S \ref{sec:Betti}.

\noindent Several examples of $C^\infty$-stable vector fields which are not (CF) are known.
Such examples can be hyperbolic (for instance, geodesic flows on compact surfaces of constant negative curvature \cite{lmm:cpt}) or parabolic (for instance, horocycles flows on compact surfaces of constant negative curvature \cite{FF1} or nilflows on nilmanifolds other than tori \cite{FF2}, \cite{FF3}).
However, there are no known examples of (CF) vector fields on manifolds other than tori. A.Katok
has proposed the following 

\begin{conjecture}(\cite{H85}, \cite{K01}, \cite{K03})
 If a closed, connected, orientable manifold $M$ admits a (CF) vector field $X$, then $M$ is diffeomorphic to a torus and $X$ is smoothly conjugate to a Diophantine vector field.
\end{conjecture}
\noindent We will refer to the above conjecture as the \emph{Katok} conjecture.

\section{Globally hypoelliptic vector fields}
\label{sec:GH}
\noindent The notion of a (CF) vector field and the related Katok conjecture were introduced an studied independently of a closely related notion introduced by Greenfield and Wallach in \cite{GW}. 

\noindent  Let $\Cal D'_n (M)$ denote the space of currents of  degree $n=\text{\rm dim}(M)$ (and dimension $0$).  We remark that currents in $\Cal D_n'(M)$ are by definition continuous linear functionals on the space $C^\infty(M)$ of smooth complex-valued functions, hence the space $\Cal D_n'(M)$ coincides with the space $\Cal D'(M)$ of distributions on $M$. Any smooth $n$-form defines 
by integration a distribution on $M$, hence the space $\Omega^n(M)$ of smooth $n$-forms on $M$ 
can be naturally identified to a subspace of  $\Cal D'_n (M)$. 

\begin{definition} \cite{GW}
A smooth vector field $X$ on $M$ is \emph{globally hypoelliptic (GH)} if  $\Cal L_X U\in 
\Omega^n (M)$ implies $U\in \Omega^n(M)$  for any $U \in \Cal D'_n (M)$\,.
\end{definition}

\noindent In \cite{GW} the authors proved several basic results on (GH) vector fields. 
The fundamental result on the dynamics of (GH) vector fileds is the following non-trivial

\begin{theorem} \label{thm:GHdyn}
(\cite{GW}, Theorem 1.1) Let $X$ be a (GH) vector field. Then the flow $\{\phi_t^X\}$
is conservative. Let $\omega$ denote the $\{\phi_t^X\}$-invariant volume form. The space 
$\Cal I_X(M)$ of $X$-invariant distributions is one-dimensional and equal to $\C \omega$. 
In particular the flow $\{\phi_t^X\}$ is uniquely ergodic and minimal (strictly ergodic).
\end{theorem}

\noindent The paper then focuses on the following conjecture, proved in a few cases:

\begin{conjecture}\cite{GW}
  If a closed, connected orientable manifold $M$ admits a
  (GH)  vector field $X$, then $M$ is diffeomorphic to a
  torus and $X$ is smoothly conjugate to a constant
  Diophantine vector field.
\end{conjecture}

\noindent We will refer to the above conjecture as the \emph{Greenfield-Wallach} conjecture. In \cite{GW} the conjecture is proved in the following cases: if $M$ has dimension $n=2$; if $M$ is of the form $G/H$, $G$ a Lie group, $H$ a co-compact closed subgroup, $X$ is the projection on $M$ of a right-invariant vector field and either $G$ is compact or $G$ is a connected, simply connected $3$-dimensional Lie group and $H=\Gamma \subset G$ is a co-compact lattice. 

\noindent 
Our proof of the conjecture in dimension $3$ is based in part on a reduction to the $3$-dimensional homogeneous case:

\begin{theorem} \label{thm:GW}  (\cite{GW}, \S 2) The Greenfield-Wallach conjecture holds if 
$M=G/\Gamma$ is a homogeneous space of a $3$-dimensional, connected and simply-connected 
Lie group and the (GH) vector field $X$ on $M$ is the projection of a right invariant vector field on $G$.
\end{theorem}

\noindent Integrable function on $M$  are naturally currents of degree $0$ (and dimension $n$), that is, continuous linear funtionals on $\Omega^n(M)$ . The space $\Cal D_0'(M)$ of currents of degree $0$ (and dimension $n$) can  be identified with the space $\Cal D'_n (M)$ of distributions but the identification depends on the choice of a volume form. Let  $I_\omega:  \Cal D'_0(M) \to \Cal D_n'(M)$ 
be the standard isomorphism defined as
$$
I_\omega: U  \to U \wedge \omega    \, , \quad  U \in \Cal D'_0 (M) \,.
$$
 If $\omega$ is $X$-invariant, that is, the Lie derivative $\Cal L_X \omega=0$, the isomorphism $I_\omega$ commutes with the operator $\Cal L_X$ on currents, in particular
$$
I_\omega (\Cal L_X U ) = \Cal L_X I_\omega(U)  \,, \quad U\in \Cal D'_0 (M)\,.
$$
  
 \smallskip
 \noindent It is an exercise to prove that every (CF) vector field is (GH). The argument is based on the following simple lemma. A smooth vector field $X$ is called \emph{ conservative} if there exists a smooth $X$-invariant volume form on $M$.
\begin{lemma} \label{lemma:GHequiv} A smooth conservative vector field $X$ on $M$ is (GH) if and only if $\Cal L_XU  \in C^{\infty}(M)$ implies $U\in C^\infty(M)$ for any current $U \in \Cal D'_0(M)$.
\end{lemma} 

\begin{proposition} Every (CF) vector field  is (GH).
\end{proposition}

\begin{proof} If $X$ is (CF) than $X$ there exists an $X$-invariant smooth volume form $\omega$ on $M$ (see Prop. \ref{prop:CFdyn}). Let $U \in \Cal D'_0(M)$ be such that $XU=f \in C^{\infty}(M)$. 
Since $X$ is (CF),  there exist a constant $c_f \in \C$ and a solution $u \in C^{\infty}(M)$ of the equation 
$Xu=f-c_f$. It follows that $X(U-u) = c_f $ in $\Cal D'_0(M)$, hence $c_f=0$. In fact, since $\Cal L_X\omega=0$, 
$$
\text{\rm vol}_\omega (M) c_f = \<c_f, \omega\> = \<X(U-u), \omega\> = \<U-u, \Cal L_X\omega\> =0\,. 
$$
However, since $X$ is (CF), the kernel of $\Cal L_X$ on $\Cal D'_0(M)$ is trivial, equal to the subspace of constant functions. It follows that  $U-u\omega \in \C$ and $U \in C^\infty(M)$.
\end{proof}

\noindent The converse statement is less evident. In terms of the notion of a $C^\infty$-stable vector field (see Definition \ref{def:stable}), the contribution of Greenfield and Wallach in this direction (see \cite{GW}, Prop. 1.5) can be formulated as follows: 
\begin{proposition} 
Every $C^\infty$-stable (GH) vector field is (CF).
\end{proposition}
\begin{proof} If $X$ is $C^\infty$-stable, then the cohomological equation $Xu=f$ has a solution
$u\in C^\infty(M)$ for every $f\in \Cal I_X (M) ^\perp \subset C^\infty(M)$. If $X$ is (GH), the space
$\Cal I_X (M)$ is one-dimensional, hence $X$ is (CF).
\end{proof}
\noindent On the basis of special cases, it natural to ``suspect that the range of a (GH) vector field is always closed'' in $\Cal D'_0(M)$ (see the Note at the end of \S 1 in \cite{GW}), that is, that every $(GH)$ vector field is $C^\infty$-stable, hence it is $(CF)$. However, this question has remained open until recently. In \cite{GW2}, Greenfield and Wallach proved a partial converse which can be formulated
as follows:
\begin{proposition} 
If  $X$ is  volume preserving, $C^\infty$-stable and if the space of $X$-invariant distributions $\Cal I_X(M) \subset \Omega^n(M)$, then $X$ is a (GH) vector field.
\end{proposition}
\begin{proof}  Let $U\in \Cal D_0'(M)$ be such that $\Cal L_X U =f  \in C^\infty(M)$.  
Since $\Cal I_X(M) \subset \Omega_n(M)$, it follows that $f\in \Cal I_X(M)^\perp$ and, since $X$ is 
$C^\infty$-stable, there exists $u\in C^\infty(M)$ such that $Xu=f$. Let $\omega$ denote the 
$X$-invariant volume form. The distribution $(U-u)\wedge \omega \in  \Cal I_X (M)\subset \Omega^n(M)$, hence $U \in C^\infty(M)$. Thus $X$ is a (GH) vector field by Lemma \ref{lemma:GHequiv}.
\end{proof}

\noindent  In 2000 Chen and Chi \cite{CC00} have published a paper based on the result that (GH) vector fields on tori are always $C^\infty$-stable. Their argument is essentially correct and generalizes word for word  to any compact manifold. We owe this remark to F. Rodriguez-Hertz.

\begin{theorem} (after Chen and Chi \cite{CC00})  Every (GH) vector field on $M$ is $C^\infty$-stable, hence it is (CF).
\end{theorem}

\begin{proof} Let $\omega$ be the (normalized) $X$-invariant volume form and let $L^2(M,\omega)$
be the standard Hilbert space of square-integrable functions with respect to the $X$-invariant
volume. Let $\{ H^s(M) \vert s\in \R\}$ be the standard family of Sobolev spaces on the compact 
manifold $M$. We remark that the space 
$$
H^\infty(M)= \cap_{s\in \N} H^s(M) = C^{\infty}(M) \,,
$$ 
endowed with the sequence of Sobolev norms $\{ \Vert \cdot \Vert_n \vert n\in \N\}$, is a 
Fr\'echet space. Let $L: H^{-1}(M) \to H^{\infty}(M)$ be the linear densely defined operator
defined as follows:
$$
L(f) = Xf   \,, \quad  f\in D(L) = H^{\infty}(M)\,.
$$
Since $X$ is (GH), the  operator $L$ on $H^{-1}(M)$ is closed on the
(dense) domain $D(L)=H^{\infty}(M)$, hence its graph $G_L=\{ (f, Lf) \vert f\in D(L)\}$ is a 
closed subspace of the Fr\'echet space $H^{-1}(M) \times H^{\infty}(M)$. The linear operator 
$\pi: G_L \to  L^2(M,\omega)$ defined as 
$$
\pi (f, Lf) =f \,, \quad f \in D(L)= H^\infty(M)\,.
$$
is closed operator on the Fr\'echet space $G_L$ (a closed subspace of the Fr\'echet space
$H^{-1}(M) \times H^{\infty}(M)$), hence it is bounded by the closed graph theorem. It follows that
there exist $s\in \N$ and  a constant $C>0$ such that, for any $f\in H^\infty(M)$,
\begin{equation}
\label{eq:mainest1}
\Vert f \Vert_0   \leq  C \left (   \Vert f  \Vert _{-1}  +   \Vert Xf \Vert_s \right) \,.
\end{equation}
We claim that there exists $C' >0$ such that the following estimate holds:
\begin{equation}
\label{eq:mainest2}
\Vert f \Vert_0   \leq  C'  \Vert Xf \Vert_s \,, \quad \text{ \rm for all }  f \in H^\infty(M) \text{ such that }
 \,\,\int_M f\omega =0 \,.
\end{equation} 
Let us assume that the claim does not hold. Then there exists a sequence $\{f_j \vert j\in\N\}$
in $H^\infty(M)$ such that $\int_M f_j\omega \equiv 0$ and 
$$
\Vert f_j \Vert_0 \equiv 1   \,, \quad  \Vert X f_j \Vert_s \to 0 \,.
$$
Since $L^2(M,\omega)$ is a (separable) Hilbert space, there exists a subsequence 
$\{f_{j_k} \vert k\in\N\}$ weakly convergent to $f\in L^2(M,\omega)$. Since  $Xf_j \to 0$ 
in $H^s(M)$, it follows that $Xf=0$ in $\Cal D'_0(M)$, hence $f\in \C$ is a constant
function. However, $\int_M f_{j_k} \omega \to  \int_M f  \omega =0$, implies that
$f=0$. By Rellich embedding theorem, the embedding $L^2(M,\omega) \to 
H^{-1}(M,\omega)$ is compact, hence $f_{j_k} \to 0$ strongly in $H^{-1}(M)$. 
Finally, by estimate \pref{eq:mainest1} we have
$$
\Vert f_{j_k}  \Vert_0   \leq  C \left (\Vert f_{j_k} \Vert _{-1}  +   \Vert Xf_{j_k} \Vert_s \right)\,,
$$
hence $f_{j_k} \to 0$ in $L^2(M, \omega)$ contradicting the assumption that $\Vert f_j \Vert_0 
= 1$ for all $j\in \N$. The claim is therefore proved. 

\smallskip
\noindent Since $L^2(M,\omega)$ is complete and  $X$ is (GH), it follows immediately from the estimate \pref{eq:mainest2} that $X$ is $C^\infty$-stable, hence it is (CF) 
\end{proof}

\smallskip
\noindent The above results can be summarized as follows:
\begin{theorem} \label{thm:equiv} 
Let $X$ be a smooth vector field on a closed connected manifold $M$. The following statements are equivalent:
\begin{enumerate}
\item $X$ is (GH);
\item $X$ is (CF);
\item $X$ is volume preserving, $C^\infty$-stable and all $X$-invariant distributions are smooth
$n$-forms.
\end{enumerate}
\end{theorem}

\noindent We became aware of the the paper \cite{CC00} and of its main result (that (GH) vector fields
on tori are smoothly conjugate to Diophantine vector fields) by reading the paper \cite{HH}. However, the question whether every (GH) vector field is $C^\infty$-stable, hence (CF) was still proposed as an open question in the paper \cite{FF3}. Only recently, F. Rodriguez-Hertz has informed us that \cite{CC00} is actually based on a proof (for the toral case) that (GH) vector fields are (CF). The authors were apparently not aware of a paper of Luz and Dos Santos  \cite{LdS} whose methods can be adapted to prove that  every (CF) vector field on a torus is smoothly conjugate to a constant  Diophantine vector field (see \cite{Ko}, \S 2.2) and give a (quite convoluted) independent proof.

\section{The first Betti number}
\label{sec:Betti}

\noindent  In this section we will outline recent progress on the Katok (and Greenfield-Wallach)
conjecture due to F. and J. Rodriguez-Hertz \cite{HH}. Their main result can be stated as follows:
a (CF) flow on a closed, connected manifold $M$ has a smooth Diophantine factor on a torus of the dimension of the first Betti number of $M$. In particular, if the first Betti number is equal to the dimension of $M$, then $M$ is diffeomorphic to a torus and the (CF) vector field is smoothly conjugated to a Diophantine vector field. 

\begin{lemma} \label{lemma:proj} (\cite{HH}, Prop. 1.3)  Let $p:M \to N$ be a smooth fibration and let 
$Y$ be a smooth vector fields on $N$ such that $Y=p_\ast (X)$. If $X$ is (CF), then $Y$ is (CF).
\end{lemma}
\begin{proof} Let $f\in C^\infty(N)$ and let $g=f\circ p$. There exists a constant $c\in \C$ and a function $v\in C^\infty(M)$ such that $Xv= g -c$.  For any $y\in N$, the set $M_y =p^{-1}\{y \}$  is a smooth submanifolds of codimension equal to the dimension of $N$.  Let $\omega$ be the $X$-invariant volume form on $M$ and let $\omega_y$ be the restriction of $\omega$ to $M_y$. The form $\omega_y$ is a volume form on $M_y$ and the function $w:M\to \R$ defined for all $x\in M$ as  
$$
w(x):= \text{ \rm vol}_{\omega_{p(x)}}(M_{p(x)})
$$
is an $X$-invariant smooth function, hence it is constant equal to $w\in \R^+$. 
Let $u\in C^\infty(N)$ be defined as
$$
u(y) := w^{-1} \int_{M_y}  v \omega_y  \,, \quad \text{\rm  for all }\, y\in N \,.
$$
A computation shows that
$$
\Cal L _Y u (y) = w^{-1}  \int_{M_y} \Cal L_X v  \omega_y = \int_{M_y} (f\circ p -c) \omega_y =
f(y) -c \,,
$$
hence $u\in C^\infty(N)$ is a solution of the cohomological equation $Yu=f-c$. It follows that $Y$ is (CF).
\end{proof}

\noindent The results of \cite{HH} are based on the following simple but crucial idea:

\begin{lemma} \label{lemma:key} Let $X$ be a (CF) vector field on $M$. There exists a continuous linear operator $j_X : \Omega^1(M)\to \Omega^1(M)$ on the space of $1$-forms with the following properties. Let $\omega$ be the normalized $X$-invariant volume form on $M$. For every $\eta \in   \Omega^1(M)$,

\begin{enumerate}
\item $\imath_X  j_X(\eta) \equiv  \int_M \imath_X \eta \,\omega $;
\item the $1$-form $j_X (\eta) -\eta$ is exact; 
\item if $\imath_X\eta$ is constant, then $j_X (\eta)$ is $X$-invariant.
\end{enumerate}
In particular, any de Rham cohomology class $c\in H^1(M,\R)$ has an $X$-invariant representative, 
that is, there exists $\eta\in \Omega^1(M)$ such that 
$$
c= [\eta] \in H^1(M,\R)   \quad \text{ \rm and } \quad \Cal L_X \eta =0 \,.
$$  
\end{lemma}
\begin{proof} Since $X$ is (CF) there exists a linear operator $u_X: \Omega^1(M) \to C^\infty(M)$ such that for every $\eta\in \Omega^1(M)$ the function $u:=u_X(\eta)\in C^\infty(M)$ is the 
unique zero-average solution of the cohomological equation  
\begin{equation}
\label{eq:CEcontr}
Xu = \imath_X\eta \,-\,  \int_M \imath_X \eta\,\omega\,.
\end{equation}
Let $j_X:\Omega^1(M)\to \Omega^1(M)$ be the operator defined as
$$
j_X(\eta) := \eta -du_X(\eta) \, , \quad \text{ \rm for every }\, \eta \in \Omega^1(M)\,.
$$
The operator $j_X$ is well-defined, linear and it is continuous since the operator $u_X$  is  
continous by the open mapping theorem, 

\smallskip
\noindent Property $(1)$ and $(2)$ hold by definition. If $ \imath_X\eta$ is constant,
$$
\Cal L_X j_X (\eta) = \Cal L_X \eta - d\Cal L_Xu_X(\eta)= \imath_X d\eta + d\imath_X \eta - d\imath_X \eta =0 \,,
$$
hence $j_X(\eta)$ is $X$-invariant. Finally, if $\eta$ is closed, the form $j_X(\eta)$ is closed 
by $(2)$ and it is $X$-invariant by $(3)$. Hence every de Rham cohomology class has an 
$X$-invariant representative.
\end{proof}

\noindent Let $\beta_1(M):= \text{ \rm dim}\,H^1(M, \R)$ be the first Betti number of the 
closed, connected,  orientable $n$-dimensional manifold $M$.

\begin{theorem}
\label{thm:HH}   Let $X$ be a (CF) vector field on $M$. There exist a smooth fibration $p:M \to \T^{\beta_1(M)}$ and a Diophantine vector field $Y$ on the torus $\T^{\beta_1(M)}$ such that $p_\ast(X)= Y$. It follows that $\beta_1(M) \leq n$ and if equality holds then $M$ is diffeomorphic to the torus $\T^n$ and $X$ is smoothly conjugate to a Diophantine vector field on $\T^n$.
\end{theorem}

\begin{proof} [outline] Let $\beta:= \beta_1(M)$ and let $\{c_1, \dots, c_\beta\} \subset H^1(M, \Z)$ be
an integer basis of the de Rham cohomology. By Lemma \ref{lemma:key} there exists a system
$\{ \eta_1, \dots , \eta_\beta\} \subset \Omega^1(M)$ of closed, $X$-invariant $1$- forms such that 
$c_k = [\eta_k] \in H^1(M, \R)$ for all $k\in \{1, \dots, \beta\}$. Let $x_0 \in M$ and let $p: M \to \T^\beta$
be the map defined as follows:
\begin{equation}
\label{eq:defp}
p(x) = \left( \int_{x_0} ^x \eta_1 , \dots, \int_{x_0}^x \eta_\beta \right)  \in \T^\beta \,, \quad \text{ \rm for
all } \, x\in M\,.
\end{equation}
The map $p$ is by definition well-defined and smooth. In fact, each of the integrals in \pref{eq:defp} are independent modulo $\Z$ on the choice of the path joining the base point $x_0$ to $x\in M$.
Since the forms ${ \eta_1, \dots , \eta_\beta} $ are closed and $X$-invariant , the vector field $Y:=p_\ast(X)$ is well-defined and constant on $p(M)$, in fact, it is equal to
$$
\left(\imath_X \eta_1, \dots,   \imath_X \eta_\beta  \right) \,.
$$
Since the flow $\{\phi_t^X\}$ is minimal and $p_\ast(X)$ is a constant vector field, it is then possible to prove: $(a)$ the range of $p$ is a closed subgroup of $\T^\beta$, hence it is a sub-torus $\T^\alpha \subset \T^\beta$; $(b)$ by Sard's theorem the map $p:M \to \T^\alpha$ has constant maximal rank, hence $p:M\to \T^\alpha$ is a fibration; $(c)$  the map $H^1(M, \R) \to H^1(p^{-1}(\{t\}), \R)$ is trivial for any $t\in \T^\alpha$, hence $\alpha=\beta_1(M)$ and the sub-torus $\T^\alpha= \T^\beta$ ($p$ is  surjective) ; $(d)$ by Lemma \ref{lemma:proj}, the constant vector field $Y=p_\ast(X)$ on $\T^\beta$ 
is (CF), hence Diophantine.
\end{proof}

\section{The case of $3$-manifolds}

\noindent In this section we prove the Katok conjecture (hence by Theorem \ref{thm:equiv}
the Greenfield-Wallach conjecture as well) by the following method. We prove by contradiction 
that if $M$ is a closed, connected orientable $3$-manifolds with first Betti number $\beta_1(M)<3$, 
then there is no (CF) vector field on $M$. The conjecture then follows from Theorem
\ref{thm:HH}. 

\noindent In case $\beta_1(M)\not =0$, we prove by an elementary argument based on Lemma 
\ref{lemma:key} and Theorem \ref{thm:HH} that any (CF) vector field has to be homogeneous.
The result of Greenfield-Wallach Theorem \ref{thm:GW} in the $3$-dimensional homogeneous 
case then implies that $M$ is a $3$-dimensional torus, a contradiction.

\noindent If $\beta_1(M)=0$, a simple key remark (which works only in dimension $3$) proves the following dichotomy: either the flow is tangent to a smooth $2$-dimensional foliations or it the Reeb flow for a smooth contact form. In the first case we again prove that the flow is homogeneous. The hard case which is left out at this point is the contact case. We can conclude the proof of the Katok conjecture by invoking the recent proof of the Weinstein conjecture by C. Taubes \cite{Taubes}. In fact, by the Weinstein conjecture every Reeb flow in dimension $3$ has at least a periodic orbit, hence cannot be uniquely ergodic, However, it seems important to develop different methods better adapted to our problem, especially in view of generalizations to higher dimensions. 

\noindent A proof of the above-mentioned results in the case $0<\beta_1(M)<3$ has been obtained independently by A. Kocsard (see \cite{Ko}, Chap. 3) with  methods  similar to those of this paper. Kocsard has also proposed an alternative proof in the case that the flow is tangent to a $2$-dimensional foliation (see \cite{Ko}, \S 4.3). His proof relies on several interesting ideas and results on the tangent
dynamics of flows on $3$-dimensional manifolds.

\subsection{ $\beta_1(M)= 2$ } \hfill
\label{dimtwo}

\smallskip
\noindent  We prove below that $M$ is a homogeneous space and the (CF) flow $\{\phi_t^X\}$ is a homogeneous flow. It then follows by the Greenfield-Wallach Theorem 
\ref{thm:GW} that $M$ is a $3$-dimensional torus, which contradicts the hypothesis on the 
dimension of the homology group.

\smallskip
\noindent By Theorem \ref{thm:HH}, there is a fibration $\pi:M \to \T^2$ such that the (CF)
vector field $X$ projects onto a constant  Diophantine vector field $\pi_*(X)$ on $\T^2$. It follows that
there exist two closed smooth $1$-forms $\eta_1$ and $\eta_2$ on $M$ such that the 
functions $\imath_X \eta_1=1$ and $\imath_X \eta_2\in \R \setminus\{0\}$  and the $2$- form
$\eta_1 \wedge \eta_2$ never vanishes on $M$. We remark that it follows that $\eta_1$ and
$\eta_2$ are invariant under the flow $\{\phi_t^X\}$, in fact the Lie derivatives
\begin{equation}
\Cal L_X \eta_i = d \imath_X \eta_i + \imath_X d\eta_i =0\,,  \quad i=1,2\,.
\end{equation}
Let $\omega$ denote the $\{\phi_t^X\}$-invariant normalized volume form on $M$. We introduce a smooth non-singular vector field $Z$ on $M$, tangent to the fibers of the fibration $\pi:M \to \T^2$ normalized
so that the following properties hold:
\begin{equation}
\label{eq:Z}
\imath_Z \eta_1 = \imath_Z \eta_2 =0  \quad \text{ \rm and } \quad  \imath_Z \omega = \eta_1\wedge \eta_2\,.
\end{equation}
This is possible since for every non-singular vector field $V$ on $M$ the kernel of the map $\imath_V:
\Omega^2(M) \to \Omega^1(M)$ is $1$-dimensional, hence it is equal to $\imath_V \Omega^3(M)$.
From properties \pref{eq:Z} it follows that $\omega$ is invariant under the flow $\{\phi^Z_t\}$, in fact
$$
\Cal L_Z \omega =d \imath_Z  \omega = d (\eta_1\wedge \eta_2) =0\,.
$$
In addition, since the $1$-forms $\eta_1$, $\eta_2$ and the volume form $\omega$ are invariant
under $\{\phi_t^X\}$ and $Z$ is uniquely determined by the properties \pref{eq:Z}, it follows that
$Z$ is $\{\phi_t^X\}$-invariant, hence the commutator $[X,Z] =0$. In fact, 
\begin{equation}
\label{eq:xz}
0=\Cal L_X  (\eta_1 \wedge \eta_2) = \Cal L_X \imath_Z \omega = \imath_{[X,Z]} \omega \,.
\end{equation}
Let $\tilde Y$ be any smooth non-singular vector field on $M$  such that 
\begin{equation}
\label{eq:tY}
\imath_{\tilde Y} \eta_1 =0 \,, \quad  \imath_{\tilde Y} \eta_2=1 \,.
\end{equation}
Since for $i=1$, $2$ 
\begin{equation}
\label{eq:xw}
0= d\eta_i (X,{\tilde Y} ) = X \eta_i({\tilde Y} ) - {\tilde Y}  \eta_i(X) - \eta_i( [X,{\tilde Y} ])\,,
\end{equation}
it follows that $ \eta_1( [X,{\tilde Y} ])= \eta_2( [X,{\tilde Y} ]) =0$, hence there exists a smooth function $f$ on $M$ such that $[X,{\tilde Y} ] =f Z$. Let $u$ be the solution of the cohomological equation
$$
Xu = f -\int_M f \omega\,.
$$
Let $Y:= {\tilde Y} - uZ$. We remark that $\imath_Y \eta_i=\imath_{\tilde Y}  \eta_i$ for $i=1$, $2$. 
We have
\begin{equation}
\label{eq:xy}
[X,Y] = [X, {\tilde Y} ] - Xu\,Z = (f- Xu)\,Z = \left( \int_M f \omega \right)\,Z\,.
\end{equation}
As in \pref{eq:xw}, the following identities hold:
\begin{equation}
\label{eq:yz}
0= d\eta_i (Y,Z) = Y \eta_i(Z) - Z \eta_i(Y) - \eta_i( [Y,Z])\,,
\end{equation}
hence there exists a smooth function $g$ on $M$ such that $[Y,Z] =g Z$. By the Jacobi identity,
by \pref{eq:xz} and \pref{eq:xy}, 
\begin{equation}
 [X,gZ] = [X, [Y,Z]] + [Z,[X,Y]] + [Y,[Z,X]] = 0
\end{equation}
hence $Xg\, Z = Xg\,Z + g[X,Z] = [X,gZ] =0$, which implies that $g$ is $\{\phi_t^X\}$-invariant, hence constant. 

\smallskip
\noindent In conclusion there exist $a$, $b\in \R$ such that 
\begin{equation}
\label{eq:commrel1}
[X,Y] = aZ  \,, \quad [Y,Z] =bZ  \quad \text{ \rm and } \quad [X,Z]=0 \,.
\end{equation}
Let $\mathfrak{g}_{a,b}$ be the (solvable) $3$-dimensional Lie algebra defined by the commutation
relation \pref{eq:commrel1} and let $G_{a,b}$ be the unique connected, simply connected Lie group
with Lie algebra $\mathfrak{g}_{a,b}$. There exists a transitive, locally free action $A:G_{a,b} \times M \to M$ of $G_{a,b}$ on $M$ by volume-preserving diffeomorphisms, defined as follows: for all $(s,t,u)\in \R^3$ and all $x\in M$,
\begin{equation}
A:\left( \exp(s X) \exp(t Y) \exp(u Z), x\right)   \to   \phi_s^X \circ \phi_t^Y \circ \phi_u^Z(x) \,,
\end{equation}
hence $M$ is a homogeneous manifold of the form $G_{a,b}/\Gamma$ for some co-compact lattice
$\Gamma$ and $\{\phi_t^X\}$ is a (CF) homogeneous flow generated by the right-invariant vector 
field $X$ on $M$. It follows by Greenfield-Wallach Theorem \ref{thm:GW} that $M$ is a
$3$-dimensional torus as claimed.

\smallskip
\noindent It is not difficult to prove that $b=0$ in the above argument, hence the group $G_{a,b}$
is nilpotent and isomorphic to the $3$-dimensional Heisenberg group. In fact, let $\eta_3$ be the smooth $1$-form on $M$ such that
\begin{equation}
\imath_X \eta_3 = \imath_Y \eta_3 =0 \quad \text{ \rm and } \quad  \imath_Z \eta_3 \equiv 1\,.
\end{equation}
As in \pref{eq:xw} and \pref{eq:yz} we compute
\begin{equation}
\label{eq:yzbis}
d\eta_3 (Z,Y) = Z \eta_3(Y) - Y \eta_3(Z) - \eta_3( [Z,Y]) \\= \eta_3( [Y,Z]) = g \,.
\end{equation}
By \pref{eq:Z} and \pref{eq:tY}, since $Y={\tilde Y} -uZ$, the following identities hold:
\begin{equation}
\imath_Y \eta_1 =\imath_Z \eta_1 =0\,,  \quad \imath_X\eta_1 = \imath_Y \eta_2=1
\quad \text{ \rm and} \quad \eta_1 \wedge \eta_2 =\imath_Z \omega\,.
\end{equation}
It follows that $ (\eta_1 \wedge d\eta_3)(X,Y,Z) = -g \, \omega(X,Y,Z)$ and, since $\eta_1$ is closed,
\begin{equation}
d(\eta_1 \wedge \eta_3)= -\eta_1 \wedge d\eta_3 =  g \,\omega\,,
\end{equation}
which implies that the constant function $g$ vanishes identically, in fact 
\begin{equation}
\int_M g \,\omega = - \int_M d(\eta_1 \wedge \eta_3) =0\,.
\end{equation}

\subsection
{ $\beta_1(M)= 1$ } \hfill
\label{dimone}

\smallskip
\noindent In this case we have a fibration $\pi: M \to \T^1$ such that $\pi_\ast (X)$ is a generator of the translations on $\T^1$. It follows that $S_\tau := \pi^{-1}(\{\tau\})\subset M$ is a smooth compact surface transverse to the flow for any $\tau\in \T^1$. Let $\Sigma_\tau$ be a connected component
of $S_\tau$ and let $f_\tau: S_\tau\to S_\tau$ be the first return map of the flow $\{\phi_t^X\}$ to the surface $\Sigma_\tau$. If $\Sigma_\tau$ is not homeomorphic to a $2$-torus $\T^2$, it can be derived from the Lefschetz fixed point theorem that $f_\tau$ must have periodic points. The argument for compact surfaces is an exercise but we refer the reader to the paper \cite{Fuller} for more general results in this vein (we are grateful to E. Pujals for this reference). It follows that the flow $\{\phi_t^X\}$ has periodic orbits, which contradicts unique ergodicity. Since $X$ is (CF), the flow  $\{\phi_t^X\}$ is
uniquely ergodic, hence $\Sigma_\tau$ is homeomorphic to $\T^2$. In this case, by the Lefschetz formula the map $f_\tau$ has no periodic points only if the linear map $(f_\tau)_*: H_1(\Sigma_\tau,\R) \to H_1(\Sigma_\tau,\R)$ has both eigenvalues equal to $1$. In this case $(f_\tau)_*$ fixes at least one (integer) homology class. The first Betti number $\beta_1(M)$ of the mapping cone $M$ of the map $f_\tau: \Sigma_\tau\to \Sigma_\tau$ is at least equal to $2$, in fact
$$
\beta_1(\hat M) = 1+ \text{ \rm dim} \text{ \rm Ker } (f_\tau^* -id ) \geq 2 \,,
$$
in contradiction with the assumption that $\beta_1(M)=1$. Hence if $\beta_1(M)=1$ there
are no (CF) vector fields on $M$.
 
 \subsection{ $\beta_1(M)= 0$ } \hfill
\label{dimzero} 

\smallskip
\noindent  If the cohomology $H^2(M,\R)=0$ and $X$ is a (CF) vector field, there exists a $1$-form $\alpha$ such that 
\begin{equation}
\label{eq:alpha1}
\imath_X \alpha \in \R  \quad \text{ \rm and } \quad  \imath_X d\alpha =0\,.
 \end{equation}
In fact, let $\eta_X:=\imath_X \omega$. Since ${\Cal L}_X \omega =0$, the form $\eta_X$ is closed.
If $H^2(M,\R)=0$, there exists a $1$-form $\theta$ on $M$ such that $d\theta=\eta_X$. Since $X$ is (CF), there exists a function $u\in C^{\infty}(M)$ such that
\begin{equation}
\imath_X \theta + Xu= \int_M \imath_X \theta \, \omega\,.
\end{equation}
The $1$-form $\alpha := \theta + du$ satisfies the required properties \pref{eq:alpha1}.

\smallskip
\noindent There are two cases:  $(a)$ $\imath_X \alpha \equiv 0$ ; $(b)$ $\imath_X \alpha \not\equiv 0$; 
In case $(a)$ it is possible to prove that $M$ is a homogeneous manifolds and $\{\phi_t^X\}$ is a
homogeneous flow, hence the Greenfield-Wallach Theorem \ref{thm:GW} implies as above that
$M$ is a $3$-torus, a contradiction. In case $(b)$ the flow generated by $X$ is the Reeb  flow for the contact structure defined by the $1$-form $\alpha$, hence it has a periodic orbit by the Weinstein conjecture, recently proved by C. Taubes \cite{Taubes}. However, every (CF) flow is 
volume preserving and uniquely ergodic, hence it cannot have periodic orbits.

\smallskip
\noindent Let us prove that in case $(a)$ $M$ is a homogeneous manifolds and $\{\phi_t^X\}$ is  a
homogeneous flow. Let $\alpha$ be a smooth $1$-form such that $d\alpha \not\equiv 0$ and
\begin{equation}
\label{eq:alpha2}
\imath_X \alpha= \imath_X d\alpha =0 \,.
\end{equation}
It follows that  $\alpha \wedge d\alpha =0$ and $\alpha$ is $\{\phi_t^X\}$-invariant, that is,
\begin{equation}
 \label{eq:LXalpha}
\Cal L_X \alpha= d \imath_X \alpha + \imath_X d\alpha = 0 \,.
\end{equation}
Since the flow $\{\phi_t^X\}$ is uniquely ergodic, it follows that the form $\alpha$ is everywhere non-singular and there exists $c\in \R\setminus \{0\}$ such that $d\alpha =c\,\eta_X$. In fact, by 
\pref{eq:alpha2} there exists $f \in C^{\infty}(M)$ such that $d\alpha= f \eta_X$. The function $f$ is $\{\phi_t^X\}$-invariant, hence  constant.  It is therefore possible to normalize $\alpha$ so that
$d\alpha = \eta_X$. By \pref{eq:alpha2} it also follows that  $\alpha \wedge d\alpha =0$,  hence
there exists a smooth $1$-form $\tilde \beta$ (not unique) such that
 \begin{equation}
 \label{eq:tbeta}
 d\alpha = {\tilde \beta} \wedge \alpha \,.
 \end{equation}
 Let us compute the Lie derivative $\Cal L_X {\tilde\beta}$. By formulas \pref{eq:LXalpha} and 
  \pref{eq:tbeta}
 \begin{equation}
 \label{eq:LXtbeta}
  0= \Cal L_X  d\alpha =  \Cal L_X {\tilde \beta} \wedge \alpha  \,,
 \end{equation} 
 hence there exists a smooth function $g$ on $M$ such that $\Cal L_X {\tilde \beta}= g \alpha$.
 Let $v$ be the solution of the equation 
 $$
 Xv + g = \int_M  g \, \omega =a \in \R 
 $$
 and let $\beta := \tilde \beta + v\alpha$. We then have
 \begin{equation}
 \label{eq:LXbeta}
  {\Cal L}_X  \beta =  {\Cal L}_X {\tilde \beta}  +  Xv \alpha = (g+Xv) \alpha 
  = a \alpha \,.
 \end{equation} 
 We remark that  
  \begin{equation}
 \label{eq:dalpha}
 d\alpha = {\tilde \beta} \wedge \alpha  = \beta \wedge \alpha\,,
 \end{equation}
 hence $0= \imath_X d\alpha = (\imath_X \beta) \alpha $ which implies 
  \begin{equation}
 \label{eq:iXbeta}
 \imath_X \beta = \imath_X \alpha = 0\,.
 \end{equation}
 
 \smallskip
 \noindent Let $\tilde \gamma$ be any smooth $1$-form such that $\imath_X {\tilde \gamma}\equiv 1$. Since 
 $$
 \omega = \imath_X \omega \wedge {\tilde \gamma}= d \alpha \wedge {\tilde \gamma} = \alpha\wedge
 \beta \wedge {\tilde \gamma} \,,
 $$
 the forms $\alpha$, $\beta$ and ${\tilde \gamma}$ are linearly independent at all $x\in M$,
 hence there exists smooth functions $h_1$, $h_2$, $h_3$ on $M$ such that
 \begin{equation}
 \label{eq:LXtgamma}
 \Cal L_X {\tilde \gamma} = h_1 \alpha + h_2\beta + h_3  {\tilde \gamma} \,.
 \end{equation}
 Since $ d \imath_X  {\tilde \gamma}= 0$, it follows that $ \Cal L_X {\tilde \gamma} = \imath_X
 d{\tilde \gamma}$, hence
 $$
 h_3= \imath_X\,  \Cal L_X {\tilde \gamma} \equiv  0 \,.
 $$
 Let $w_1$ and $w_2$ be the smooth solutions of the cohomological equations
 \begin{equation}
 \label{eq:w1w2}
 \begin{aligned} 
 Xw_1 + h_1 + a w_2 &=  \int_M h_1 \omega + a\int_M w_2\,\omega  = b \in \R  \\
  Xw_2 + h_2 &= \int_M h_2 \, \omega =c \in \R \,.
 \end{aligned}
 \end{equation}
 The above equations can be solved since $X$ is a (CF) vector field. In fact, the
 second equation does not depend on the first equation, hence it has a solution $w_2 
 \in C^\infty(M)$. Once the solution $w_2$ has been chosen, the first equation becomes a
 cohomological equation for $w_1$ and can also be solved.  
 \smallskip
 \noindent Let $\gamma := {\tilde \gamma} + w_1 \alpha + w_2 \beta$. We remark that $\imath_X
 \gamma \equiv \imath_X{\tilde \gamma} \equiv 1$. A computation yields
 \begin{equation}
 \begin{aligned}
 {\Cal L}_X \gamma &= {\Cal L}_X {\tilde \gamma} + (Xw_1) \,\alpha + (Xw_2)\, \beta + w_2 
 \,{\Cal L} _X\beta \\ 
 &= (Xw_1 + h_1+ a w_2) \,\alpha + (Xw_2 +h_2) \,\beta= b \,\alpha + c \,\beta \,.
 \end{aligned}
 \end{equation}
 We have thus proved that there exists $a$, $b$, $c\in \R$ such that
 \begin{equation}
 \label{eq:LXids}
 {\Cal L}_X \alpha = 0 \,, \quad {\Cal L}_X\beta =a\,\alpha \quad \text{ \rm and }
 \quad {\Cal L}_X \gamma = b\,\alpha + c\,\beta\,.
 \end{equation}
 The above identities show that the flow $\{\phi_t^X\}$ is `homogeneous'. In order to
 prove that the manifold $M$ has an homogeneous structure, we will compute the
 differentials of the forms $\alpha$, $\beta$ and $\gamma$. 
 Since 
 $$
 \imath_X \alpha =\imath_X \beta =0 \quad \text{ \rm and } \quad  \imath_X \gamma =1 \,,
 $$
 it follows from \pref{eq:LXids} that $\imath_X d\beta = \Cal L_X \beta$ and 
 $ \imath_X d\gamma = \Cal L_X \gamma$, hence there exist smooth functions
 $r_1$, $r_2\in C^\infty(M)$ such that
 \begin{equation}
 \label{eq:dbetagamma}
 \begin{aligned}
 d\beta &= -a(\alpha \wedge \gamma) +  r_1 (\alpha\wedge\beta) \,, \\
 d\gamma &= -b(\alpha \wedge \gamma) - c (\beta\wedge \gamma) + 
 r_2 (\alpha\wedge\beta)\,.
 \end{aligned}
 \end{equation}
  Since the forms $\alpha$ and $d\alpha =\beta \wedge \alpha$ are  $\{\phi_t^X\}$-invariant, a computation yields:
  \begin{equation}
  \label{eq:LXdbeta}
  \begin{aligned}
   \Cal L_X  d\beta &= d \imath_X d\beta + \imath_X d^2\beta=  d \imath_X d\beta = d (a \alpha )= - a \,\alpha \wedge \beta \,;\\
 \Cal L_X  d\beta &= -a \,\alpha \wedge \Cal L_X\gamma + (Xr_1)\, \alpha\wedge \beta
 = (Xr_1 - ac) \, \alpha \wedge  \beta  \,.
  \end{aligned}
 \end{equation} 
 It follows that 
 $$
 Xr_1= ac-a \in \R\,,
 $$
 which implies that $ac-a=0$, hence $a=0$ or $c=1$, and $r_1\in \R$ is a constant function. Similarly, 
 we compute
   \begin{equation}
   \label{eq:LXdgamma1}
  \begin{aligned}
   \Cal L_X  d\gamma &= d \imath_X d\gamma + \imath_X d^2\gamma=  d \imath_X d\gamma
 \\ &= b\,d\alpha + c \,d\beta =   (c\,r_1-b) \, \alpha \wedge \beta - ac \,\alpha\wedge \gamma
    \end{aligned}
 \end{equation} 
 and 
   \begin{equation}
   \label{eq:LXdgamma2}
  \begin{aligned}
   \Cal L_X  d\gamma &= -b \,\alpha \wedge \Cal L_X \gamma -c\, \Cal L_X \beta \wedge \gamma 
   -c \,\beta \wedge \Cal L_X \gamma + (Xr_2) \, \alpha \wedge  \beta\\ &= -bc\, \alpha\wedge \beta 
   -ac \,\alpha\wedge \gamma  -bc\,  \beta \wedge \alpha +  (Xr_2) \, \alpha \wedge  \beta  \\
&=   (Xr_2) \, \alpha \wedge  \beta -ac \,\alpha\wedge \gamma  \,.
    \end{aligned}
 \end{equation} 
 It follows that 
 $$
 Xr_2 = c\,r_1 -b \in \R  \,,
 $$
 which implies that $c\,r_1-b=0$ and $r_2\in \R$ is a constant function.
 
 \smallskip
 \noindent  We remark that it is possible to distinguish two cases: $(i)$ $a\not =0$ and $(ii)$ $a=0$. In case $(i)$ we can assume that $b=0$. In fact, we let 
 $$
 \gamma':= \gamma  -  \frac{b}{a} \beta \,.
 $$
 We remark that we have $\imath_X \gamma'= \imath_X\gamma \equiv 1$. We compute
 $$
 \Cal L_X \gamma' = \Cal L_X \gamma- \frac{b}{a}\, \Cal L_X \beta = c\,\beta \,.
 $$
It follows that  in case $(i)$ we can take 
 $$ 
 a\not =0\,, \quad b=0\,, \quad \text{ hence } \quad c=1\,, \quad  r_1=0 \,.
 $$
 Let us introduce the unique frame $\{X,Y,Z\}$ of the tangent bundle defined by the conditions
 \begin{equation}
  \begin{aligned}
   \imath_X \gamma &=1\,,  \quad  \text{ \rm and } \quad  \imath_X \alpha =  \imath_X \beta =0\,; \\
 \imath_Y \beta &=1\,,  \quad  \text{ \rm and } \quad  \imath_Y \alpha =  \imath_Y \gamma =0\,; \\
  \imath_Z \alpha &=1\,,  \quad  \text{ \rm and } \quad  \imath_Z \beta =  \imath_Z\gamma =0\,.
  \end{aligned}
 \end{equation}
 A computation based on the equations  \pref{eq:dalpha} and \pref{eq:dbetagamma} shows that the frame $\{X,Y,Z\}$ generates the $3$-dimensional Lie algebra characterized by the following commutation relations:
\begin{equation}
 \begin{aligned}
  \left[ X,Y \right] &= c X  \,, \\
  \left[X,Z\right] &= -b X - aY \,, \\  
  \left[Y,Z\right] &= r_2X + r_1 Y - Z \,.
  \end{aligned}
 \end{equation} 
 As in \S \ref{dimone} we conclude that $M$ is a homogeneous manifold and $\{\phi_t^X\}$
 is a (CF) homogeneous flow. By Greenfield-Wallach Theorem \ref{thm:GW} it folllows that $M$ is a $3$-dimensional torus (and $X$ is a Diophantine vector field).

\medskip
\noindent The proof of the Greenfield-Wallach and Katok conjectures in dimension $3$ is thus reduced to the proof of the Weinstein conjecture, recently announced by C. Taubes \cite{Taubes}. However, it is important in our opinion to find an alternative proof in the contact case.

\section{Some open questions and a conjecture}

\noindent The Greenfield-Wallach and Katok conjectures remain open in dimension higher than $3$
and  there are no results available in the general case other than \cite{HH}. We would like to propose
some partial problems which we think may be relevant partial steps toward a solution. The selection of such problems is quite obviously influenced by the approach we carried out in the
$3$-dimensional case. It is entirely possible that completely different ideas are needed. 

\begin{problem} Find an alternative proof that there are no (CF) contact vector fields on 
$3$-dimensional manifolds (or rational homology spheres).
\end{problem}

\begin{problem} \label{pro2} Prove the Katok conjecture for homogeneous flows \emph{in arbitrary dimensions }(that is, for homogeneous flows on closed, connected, homogeneous manifolds $M=G/\Gamma$ where $G$ is an arbitrary connected, simply connected Lie group and $\Gamma$ is a co-compact lattice).
\end{problem}

\noindent If $M$ is a nilmanifold,  that is, when $G$ is a nilpotent Lie group, Problem \ref{pro2} has been solved by the author in collaboration with L. Flaminio \cite{FF3}.

\smallskip
\noindent Finally, we remark that all known examples of volume preserving (uniquely ergodic) vector fields which fail to be (CF) have large spaces of invariant distributions with the exception of constant Liouville vector fields on tori. This suggests that the only source of lack of stabilty comes from the Liouville phenomenon on a toral factor. In particular we propose the following:

\begin{conjecture} If a closed, connected, orientable manifold $M$ admits a volume-preserving vector field $X$ such that the space $\Cal I_X(M)$ of all $X$-invariant distributions is $1$-dimensional, then $M$ is diffeomorphic to a torus. 
\end{conjecture}

\noindent Obviously, the above conjecture is stronger that the Katok conjecture. It is true in dimension $2$ and for homogenous flows in dimension $3$. It is also true for nilflows in all dimensions \cite{FF3}.
It is open in all other cases and seems to be beyond reach at the moment even in dimension $3$.

 \vskip 1cm

 \vskip 1cm 
 
 \bibliography{biblio}

\bibliographystyle{amsalpha}

\end{document}